\documentclass[11pt]{amsart}

\usepackage{amsmath,amssymb,amsthm,amsrefs,a4wide}
\usepackage{latexsym}
\usepackage{indentfirst}
\usepackage{graphicx}
\usepackage{placeins}
\usepackage{booktabs}
\usepackage{algorithm}
\usepackage{algorithmic}
\usepackage{multirow}

\theoremstyle{plain}

\theoremstyle{definition}

\theoremstyle{remark}

\renewcommand{\d}{\mathrm{d}}

\newcommand{\eps}{\epsilon}

\newcommand{\bbm}{\begin{bmatrix}}
\newcommand{\ebm}{\end{bmatrix}}
\newcommand{\R}{\mathbb{R}}
\newcommand{\C}{\mathbb{C}}
\newcommand{\D}{\mathbb{D}}

\newcommand{\p}{\partial}
\newcommand{\G}{\mathcal{G}}

\newcommand{\z}{\xi}
\renewcommand{\t}{\tau}

\newcommand{\hG}{\hat{G}}

\newcommand{\tG}{\tilde{G}}
\newcommand{\tH}{\tilde{H}}

\newcommand{\dmax}{{d_{\text{max}}}}

\newcommand{\dt}{\mathrm{d}t}
\newcommand{\dtot}{\frac{\dt}{t}}

\newcommand{\dth}{\mathrm{d}\theta}

\renewcommand{\Im}{\mathrm{Im}}

\begin{document}

\title[Analytic continuation from limited noisy Matsubara data]{Analytic continuation from limited noisy Matsubara data}

\author[]{Lexing Ying}

\address[Lexing Ying]{Department of Mathematics, Stanford University, Stanford, CA 94305}

\email{lexing@stanford.edu}

\thanks{The author thanks Lin Lin and Anil Damle for discussions on this topic.}

\keywords{Rational approximation, Prony's method, analytic continuation.}
\subjclass[2010]{30B40, 93B55.}

\begin{abstract}
  This note proposes a new algorithm for estimating spectral function from limited noisy Matsubara
  data.  We consider both the molecule and condensed matter cases. In each case, the algorithm
  constructs an interpolant of the Matsubara data and uses conformal mapping and Prony's method to
  estimate the spectral function. Numerical results are provided to demonstrate the performance of
  the algorithm.
\end{abstract}

\maketitle

\section{Introduction}\label{sec:intro}


For any non-negative spectral distribution (or function) $A(x)$ defined on $\R$, the Green's
function $G(z)$ for $z\in\C$ is given by
\[
G(z) = \frac{1}{2\pi} \int_\R \frac{1}{z-x} A(x) \d x
\]
and $A(x) = -2 \Im G(x+i 0^+)$.  For a fixed inverse temperature $\beta$, the Matsubara grid is
defined as $\{z_n=\frac{(2n-1)\pi}{\beta} i\}$. One analytic continuation problem is to recover
$A(x)$ given $G(z_n)$ at the Matsubara grid. This is known to be a highly ill-posed inverse problem
\cite{trefethen2020quantifying}. In practice, the situation is even worse due to the following
reasons. First, often only the values of $G(z)$ at a {\em limited} number of Matsubara points
$z_1,\ldots,z_N$ are provided. Second, the values of $\{G(z_n)\}$ almost always come with
noise. These two constraints make the problem even more challenging. In most quantum mechanics
computations, there are two typical cases.
\begin{itemize}
\item The molecule case, where $A(x)$ is a sum of a small number of Dirac deltas on $\R$ with support
  bounded away from zero.
\item The condensed matter case, where $A(x)$ is a positive continuous function on $\R$.
\end{itemize}

\subsection{Related work}
In computational physics and chemistry, the Matsubara Green's function data can be obtained from
finite temperature simulations using for example GW theory or quantum Monte Carlo. The spectral
function $A(x)$ describes the single-particle excitation spectrum \cite{bruus2004many}. Many methods
have been proposed for this analytic continuation problem, including Pade approximation
\cite{vidberg1977solving,beach2000reliable,schott2016analytic}, maximum entropy methods
\cite{jarrell1996bayesian,beach2004identifying,levy2017implementation,kraberger2017maximum,rumetshofer2019bayesian},
stochastic analytic continuation
\cite{sandvik1998stochastic,vafayi2007analytical,goulko2017numerical,krivenko2019triqs}, and the
more recent development based on Nevanlinna functions \cite{fei2021nevanlinna,fei2021analytical}.

This problem is also highly related to a couple of other well-studied problems, including rational
function approximation and interpolation
\cite{berljafa2017rkfit,berrut2004barycentric,beylkin2009nonlinear,gonnet2013robust,gustavsen1999rational,nakatsukasa2018aaa,wilber2021data},
approximation with exponential sums \cite{beylkin2005approximation,potts2013parameter}, and
hybridization fitting \cite{mejuto2020efficient}.

\subsection{Contributions}
This note proposes a new algorithm for both the molecule and condensed matter cases. Since the
problem is ill-conditioned, some regularization or prior information is needed.

In the molecule case, the fact that $A(x)$ is a sum of positive Dirac deltas provides a strong
prior.  Since $A(x)$ is supported away from the origin, $G(z)$ is analytic in a neighborhood of the
interval $[-\frac{(2N-1)\pi}{\beta}i,\frac{(2N-1)\pi}{\beta}i]$. The proposed algorithm proceeds by
(1) constructing an accurate interpolant of $G(z)$ over this interval, (2) using conformal mapping
to unzip this interval into a circle, and (3) applying Prony's method to identify the poles in the
unzipped domain.

In the condensed matter case, $G(z)$ is analytic in the upper half plane. The concept of
quasi-particle refers to the poles of the analytic continuation of $G(z)$ in the negative half
plane. One physically meaningful prior is that $A(x)$ can be well-approximated by a small number of
quasi-particles. Under this prior, the algorithm proceeds by (1) constructing an accurate
interpolant of $G(z)$ over the interval $[\frac{\pi}{\beta}i,\frac{(2N-1)\pi}{\beta}i]$, (2) using
conformal mapping to unzip the interval into a circle, (3) applying Prony's method to identify the
poles of the quasi-particles, and (4) finally evaluating the spectral function $A(x)$.

The rest of the note is organized as follows. Section \ref{sec:mol} describes the algorithm for the
molecule case and presents some numerical results. Section \ref{sec:cdm} is concerned with the
condensed matter case.

\section{Molecule case}\label{sec:mol}

\subsection{Algorithm}

In the molecule case, $A(x)$ and $G(z)$ take the form
\[
A(x) = \sum_j A_j \delta_{\z_j}(x), \quad G(z) = \frac{1}{2\pi} \sum_j \frac{1}{z-\z_j} A_j,
\]
where the sum is over the finite discrete support of $A(x)$. Since $A(x)$ is supported away from the
origin, let $[-\eps,\eps]$ be the maximum interval in which $A(x)$ vanishes. Due to this gap, $G(z)$
is analytic in a neighborhood of the interval $[-bi,bi]$, where $b=\frac{(2N-1)\pi}{\beta}$.

The first step is to construct an accurate interpolant of $G(z)$ for $z\in[-bi,bi]$. The interpolant
adopted here is of the form
\[
G(z) \approx \frac{1}{2\pi} \sum_k \frac{1}{z-x_k} X_k,
\]
where $\{x_k\}$ is a set of points in $(-\infty,-\eps] \bigcup [\eps,\infty)$. To motivate our
    choice of $\{x_k\}$, consider the map $z \leftrightarrow \eps/z$ that sends $(-\infty,-\eps]
  \bigcup [\eps,\infty)$ to the interval $[-1,1]$. From the classical results of approximation
    theory, the Chebyshev grid on $[-1,1]$ is a near-optimal choice for interpolation and this
    motivates the following choice for $\{x_k\}$: let $N_I=O(N)$ be an even number and define
\[
x_k = \frac{\eps}{\cos\left(\frac{k\pi}{N_I-1}\right)},\quad k=0,\ldots,N_I-1.
\]
Given $\{x_k\}$, we solve for the weights $\{X_k\}$ from the least square problem
\begin{equation}
  G(z_n) \approx \frac{1}{2\pi} \sum_k \frac{1}{z_n-x_k} X_k, \quad i=1,\ldots,N
  \label{eq:molls}
\end{equation}
Given $\{x_n\}$ and $\{X_n\}$, the interpolant denoted by $\tG(z)$ is given by
\begin{equation}
  \tG(z) \equiv \frac{1}{2\pi} \sum_k \frac{1}{z-x_k} X_k
  \label{eq:molint}
\end{equation}
is an accurate approximation to $G(z)$ on $[-bi,bi]$ that allows for {\em arbitrary sampling}.

\begin{figure}[h!]
  \begin{tabular}{ccc}
    \includegraphics[scale=0.28]{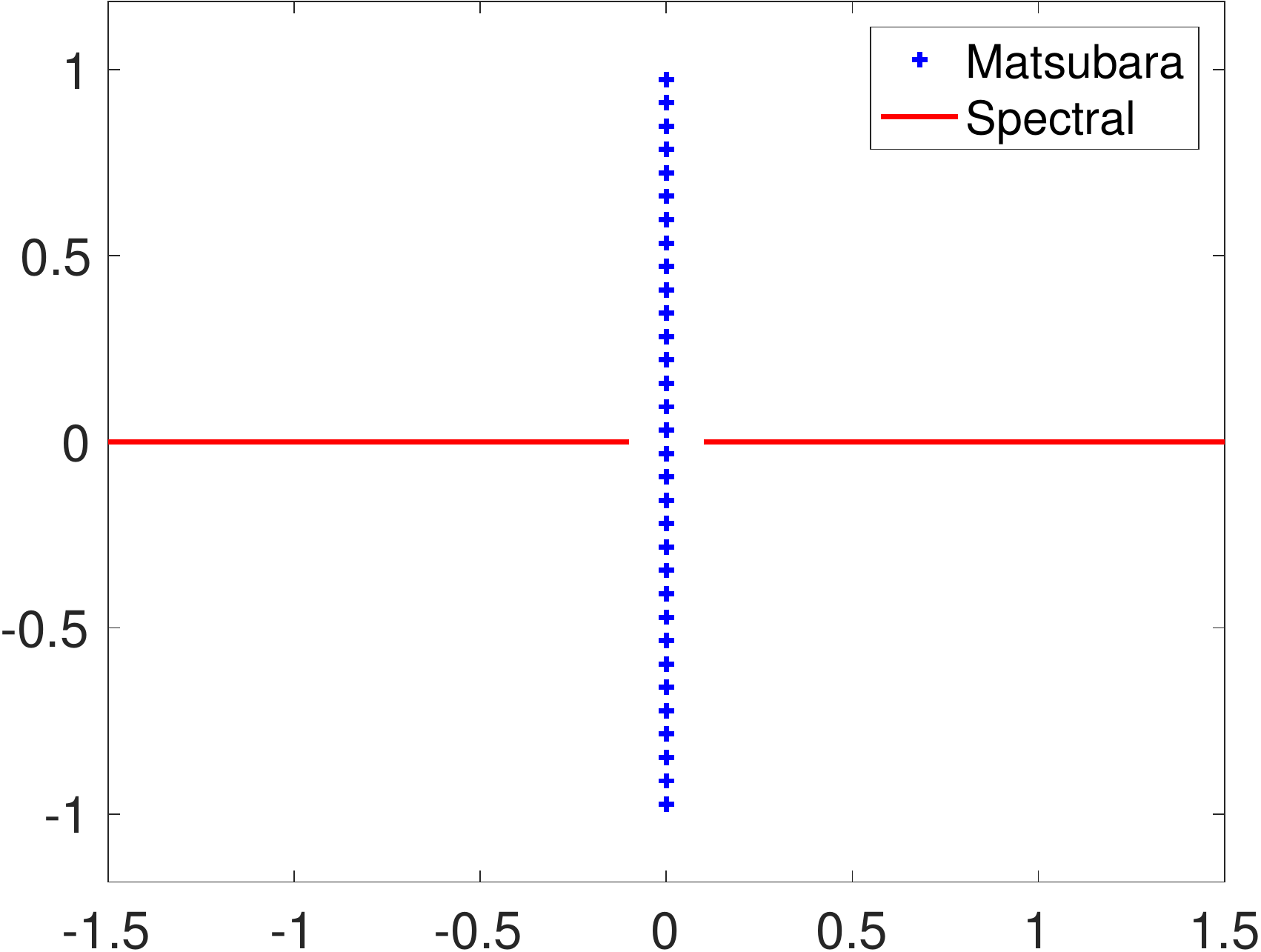}&
    \includegraphics[scale=0.28]{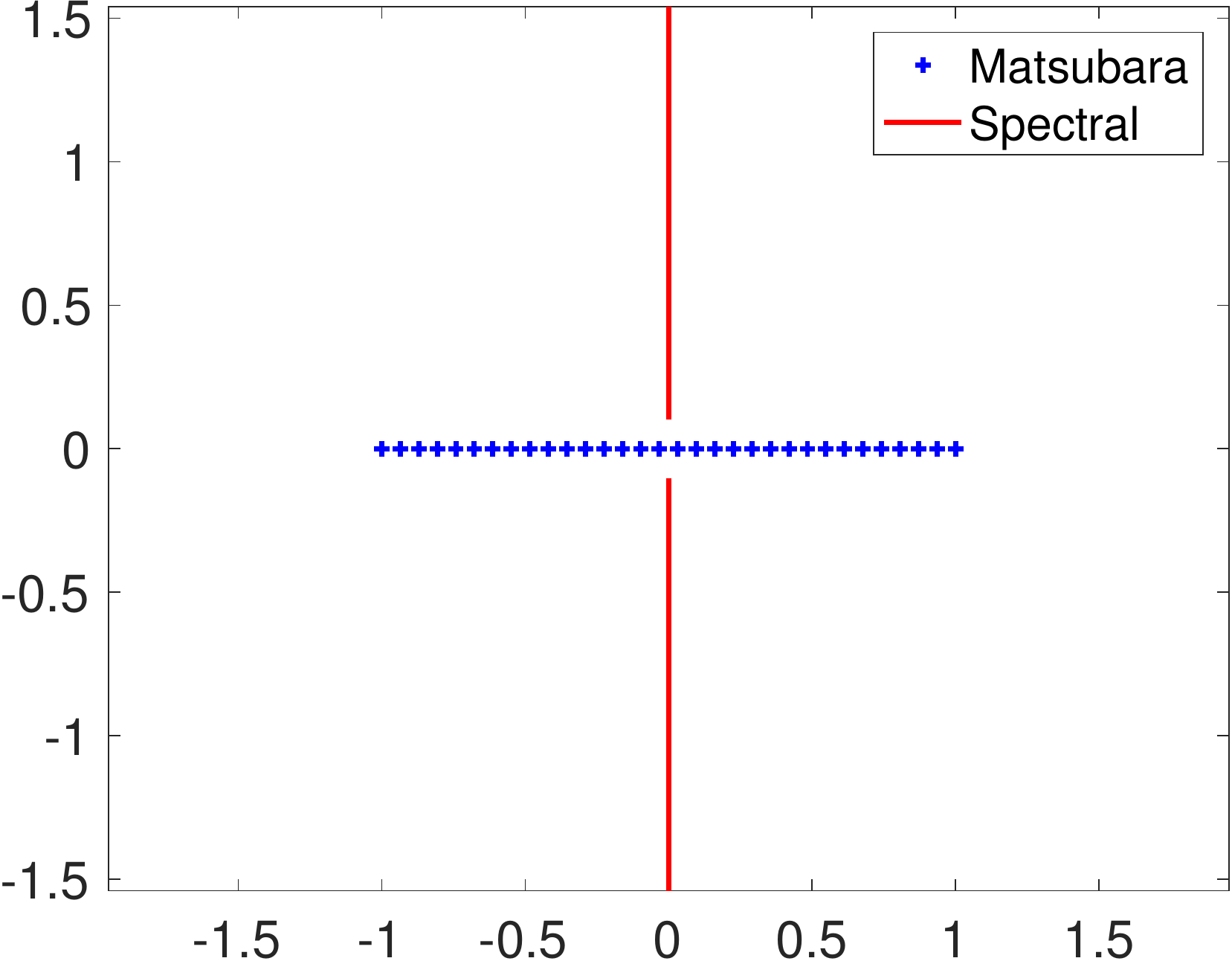}&
    \includegraphics[scale=0.28]{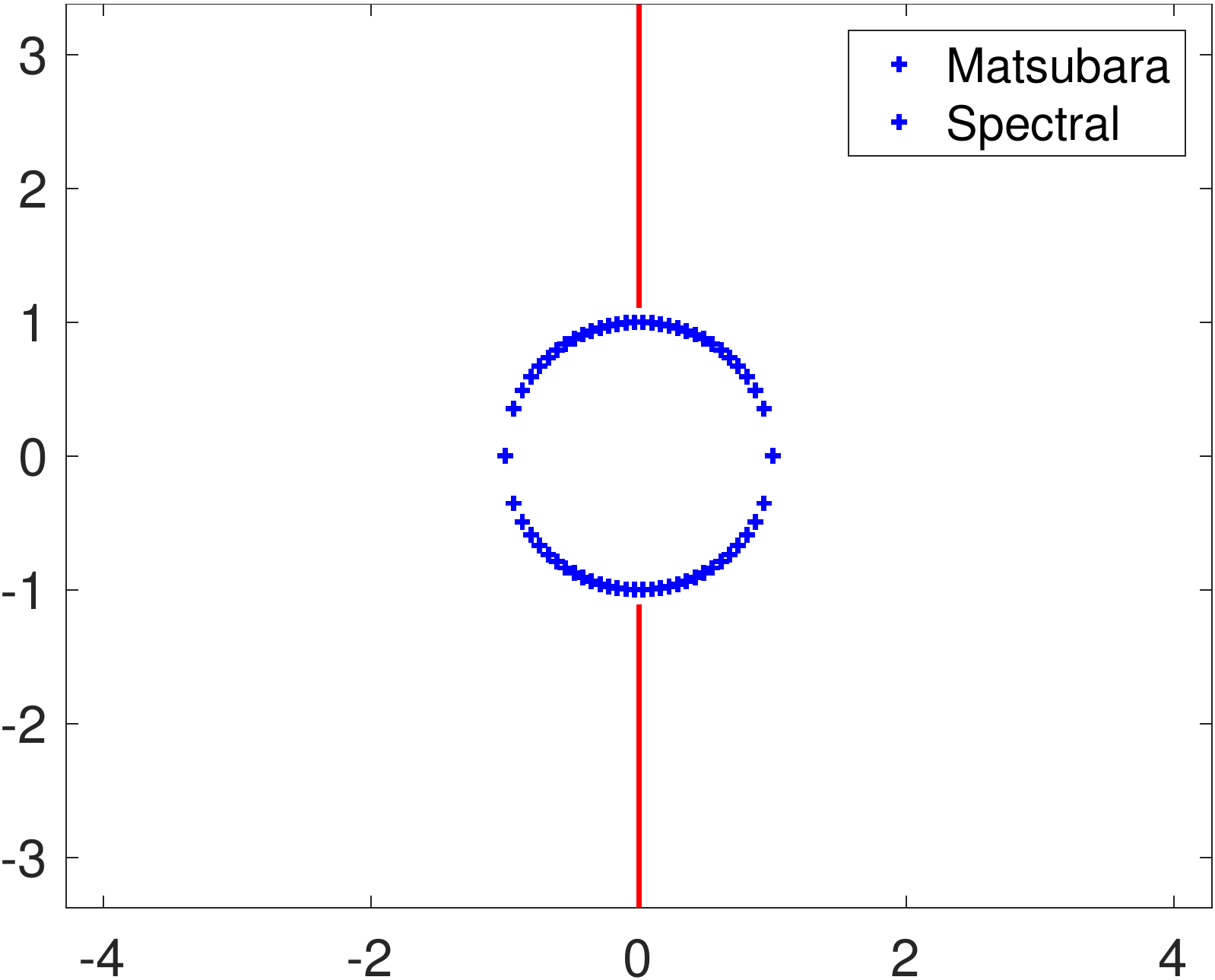}\\
    $z$ & $w$ & $t$
  \end{tabular}
  \caption{Conformal maps from $z$ to $w$ and to $t$, unzipping the interval $[-bi,bi]$ into a unit
    disk.}
  \label{fig:molunzip}
\end{figure}

The second step is to find the support $\{\z_j\}$ of the spectral distribution $A(x)$. Introduce the
following sequence of conformal mappings from $z$ to $w$ and to $t$ that unzip the interval
$[-bi,bi]$ in the $z$ plane to the unit disk $\D$ in the $t$ plane (see Figure \ref{fig:molunzip})
\begin{equation}
  w = z/(ib),\quad t = w + \sqrt{w^2 - 1}.
  \label{eq:molmap}
\end{equation}
The inverse maps are
\begin{equation}
  z = (ib) w, \quad w = \frac{1}{2} \left(t + \frac{1}{t}\right).
  \label{eq:molmapinv}
\end{equation}
In the $t$ plane, the function $G(t)\equiv G(z(t))$ is analytic outside $\D$ and takes the form
\[
G(t) = \sum_j \frac{T_j}{t-\t_j} + f(t),
\]
where $\{\t_j\}$ are the poles outside the unit disc $\D$ and $f(t)$ is analytic outside $\D$. Since
the poles $\{\z_j\}$ of $G(z)$ are mapped to the poles of $G(t)$ outside $\D$, it is sufficient to
find the poles $\{\t_j\}$ in the $t$ plane outside $\D$.

Let us consider the integrals
\begin{equation}
  \frac{1}{2\pi i}\int_{\gamma} \frac{G(t)}{t^k} \dtot
  \label{eq:integral}
\end{equation}
for integer values of $k\ge 1$, where $\gamma$ is the unit circle in the counterclockwise orientation.
For any $k\ge 1$,
\[
\frac{1}{2\pi i}\int_{\gamma} \frac{G(t)}{t^k} \dtot  = - \sum_j \frac{1}{2\pi i}\int_{\gamma_j} \frac{G(t)}{t^k} \dtot =
- \sum_j T_j \t_j^{-(k+1)},
\]
where each $\gamma_j$ is an infinitesimal circle around $\t_j$. Here, the first equality comes from
the facts that $G(t)/t^{k+1}$ is analytic outside $\gamma$ and $\{\gamma_j\}$, and $G(t)/t^{k+1}$
decays rapidly at infinity. The second infinity is because the residue of $G(t)/t^{k+1}$ at $\t_j$
is $T_j\t_j^{-(k+1)}$. This demonstrates that the integrals \eqref{eq:integral} for $k\ge 1$ contain
information about the poles outside $\D$.

Since the integral \eqref{eq:integral} is over the unit circle, it is closely related to the Fourier
transform of $G(\theta) \equiv G(e^{i\theta})$:
\begin{equation}
  \frac{1}{2\pi i}\int_{\p\D} \frac{G(t)}{t^k} \dtot = \frac{1}{2\pi i}\int_0^{2\pi} G(\theta) e^{-ik\theta}  i \dth
  = \frac{1}{2\pi} \int_0^{2\pi} G(\theta) e^{-ik\theta} \dth = \hG_k.
  \label{eq:Fourier}
\end{equation}
We emphasize that the computation of $\hG_k$ is possible because the interpolant $\tG(z)$ allows for
sampling at any $z\in[-bi,bi]$.

To recover the poles outside $\D$, we apply Prony's method to the Fourier coefficients, following
\cite{ying2022pole}. Define the semi-infinite vector
\[
\hG_+  \equiv
\begin{bmatrix}
  \hG_{1}\\
  \hG_{2}\\
  \vdots
\end{bmatrix}
\equiv
\frac{1}{2\pi i}\int_{\p\D} G(t) 
\begin{bmatrix}
  t^{-2}\\
  t^{-3}\\
  \vdots
\end{bmatrix}
\dt
\equiv
\begin{bmatrix}
  -\sum_{|\t_j|>1} T_j \t_j^{-2}\\
  -\sum_{|\t_j|>1} T_j \t_j^{-3}\\
  \vdots
\end{bmatrix}.
\]
Let us define $S$ to be the shift operator that shifts the semi-infinite vector upward (and drops
the first element). For any $\t_j$ with $|\t_j|>1$
\[
S
\begin{bmatrix}
  \t_j^{-2}\\
  \t_j^{-3}\\
  \vdots
\end{bmatrix}
=
\begin{bmatrix}
  \t_j^{-3}\\
  \t_j^{-4}\\
  \vdots
\end{bmatrix},
\quad\text{i.e.,}\quad
(S-\t_j^{-1})
\begin{bmatrix}
  \t_j^{-2}\\
  \t_j^{-3}\\
  \vdots
\end{bmatrix}
= 0.
\]
Since the operators $S-\t_j^{-1}$ all commute, 
\begin{equation}
  \prod_{\ell} (S- \t_\ell^{-1})
  \begin{bmatrix}
    \t_j^{-2}\\
    \t_j^{-3}\\
    \vdots
  \end{bmatrix}
  = 0.
\label{eq:prodL}
\end{equation}
Since $\hG_+$ is a linear combination of such semi-infinite vectors, 
\[
\prod_{\ell} (S- \t_\ell^{-1}) \hG_+ = 0.
\]
Suppose that the polynomial $\prod_{\ell} (t-\t_\ell^{-1}) = p_0 t^0 + \cdots + p_d t^d$, where the
degree $d$ is equal to the number of poles outside $\D$. Then \eqref{eq:prodL} becomes
\begin{equation}
  p_0 (S^0 \hG_+) + \cdots + p_d (S^d \hG_+) = 0,
  \quad\text{i.e.,}\quad
\begin{bmatrix}
  \hG_{1} & \hG_{2} & \cdots & \hG_{d+1}  \\
  \hG_{2} & \hG_{3} & \cdots & \hG_{d+2} \\
  \vdots & \vdots & \vdots & \vdots
\end{bmatrix}
\begin{bmatrix}
  p_0\\
  \ldots\\
  p_d
\end{bmatrix}
= 0.
\label{eq:lsL}
\end{equation}
This implies that the number of poles outside $\D$ is equal to the smallest value $d$ such that the
matrix in \eqref{eq:lsL} is rank deficient. $(p_0,\ldots,p_d)$ can be computed as a non-zero vector
in the null-space of this matrix and the roots of
\[
p(t) = p_0 t^0 + \ldots p_d t^d
\]
are $\{\t_j^{-1}\}$. Taking inverse of these roots gives the poles $\{\t_j\}$ outside $\D$. Applying
the inverse maps \eqref{eq:molmapinv} from $t$ to $w$ and to $z$ leads to the poles $\{\z_j\}$,
i.e., the support of $A(x)$ in the $z$ plane.

In the third step, with the poles located we solve the constrained optimization problem
\begin{equation}
  \min_{A_j\ge 0} \sum_i \left|G(z_i) - \frac{1}{2\pi} \sum_j \frac{A_j}{z_i-\z_j} \right|^2
  \label{eq:molopt}
\end{equation}
for $\{A_j\}$.

To implement this algorithm, we need to take care several numerical issues.
\begin{itemize}
  
\item Computation of the weights $X_k$ in \eqref{eq:molls} requires least square solution. This is
  done by a pseudo-inverse with relative singular value cutoff at $10^{-8}$.
  
\item The semi-infinite matrix in \eqref{eq:lsL}. In the implementation, pick a value $\dmax$
  that is believed to be the upper bound of the number of poles and form the matrix
  \begin{equation}
    H =
    \begin{bmatrix}
      \hG_{1} & \hG_2 & \cdots & \hG_{\dmax} \\
      \hG_{2} & \hG_3 & \cdots & \hG_{(\dmax+1)} \\
      \vdots & \vdots & \vdots & \vdots\\
      \hG_{l} & \hG_{l+1} & \cdots & \hG_{(\dmax+l-1)}
    \end{bmatrix}
    \label{eq:H}
  \end{equation}
  with $l$ satisfying $l \ge \dmax$. In practice, $l=\dmax$ is enough.

\item Numerical estimation of the rank $d$ in \eqref{eq:lsL}. To address this, let
  $s_1,s_2,\ldots,s_{\dmax}$ be the singular values of the matrix $H$. The numerical rank is chosen
  to be the smallest $d$ such that $s_{d+1}/s_1$ is below the noise level.
  
\item Computation of the vector $p$. We first compute the singular value decomposition (SVD) of
  \[
  \begin{bmatrix}
    \hG_{1} & \hG_2 & \cdots & \hG_{d+1} \\
    \hG_{2} & \hG_3 & \cdots & \hG_{d+2} \\
    \vdots & \vdots & \vdots & \vdots\\
    \hG_{l} & \hG_{l+1} & \cdots & \hG_{d+l}
  \end{bmatrix},
  \]
  respectively for \eqref{eq:lsL}. $p$ is then chosen to be the last column of the
  $V$ matrix.
  
\item The matrix $H$ in \eqref{eq:H} requires the Fourier transform $\hG_k$ from $k=-(\dmax+l-1)$ to
  $(\dmax+l-1)$. Computing the integrals in \eqref{eq:Fourier} requires evaluating $G(\theta)$ at
  quadrature points. With the interpolant
  \[
  \tG(\theta)\equiv \tG(z(\exp(i\theta)))
  \]
  available, we can sample at any point $\theta\in[0,2\pi]$.  Choose an even $N_s \ge 2(\dmax+l)$
  and define $\theta_n = \frac{2\pi n}{N_s}$ for $n=0,\ldots,N_s-1$. Using samples
  $\{\tG(\theta_n)\equiv \tG(z(\exp(i\theta_n)))\}$ at the points $\{\theta_n\}$ corresponds to
  approximating the integrals in \eqref{eq:Fourier} with the trapezoidal rule. The trapezoidal rule
  is exponentially convergent for smooth functions when the step size $h = \frac{2\pi}{N_s}$ is
  sufficient small. In the current setting, this corresponds to
  \[
  h \ll \frac{\eps}{b}, \quad\text{i.e.,}\quad  N_s \gg \frac{b}{\eps}.
  \]
  Applying the fast Fourier transform to $\{\tG(\theta_n)\}$ gives accurate approximations to
  $\{\hG_k\}$ for $k=-\frac{N_s}{2},\ldots,\frac{N_s}{2}-1$. Among them,
  $\hG_{-(\dmax+l-1)},\ldots,\hG_{(\dmax+l-1)}$ are used to form the $H$ matrix in \eqref{eq:H}.

\item The constrained optimization problem \eqref{eq:molopt} is solved with CVX \cite{cvx}.
  
\end{itemize}


\subsection{Numerical results}

\begin{figure}[h!]
  \begin{tabular}{cc}
    \includegraphics[scale=0.3]{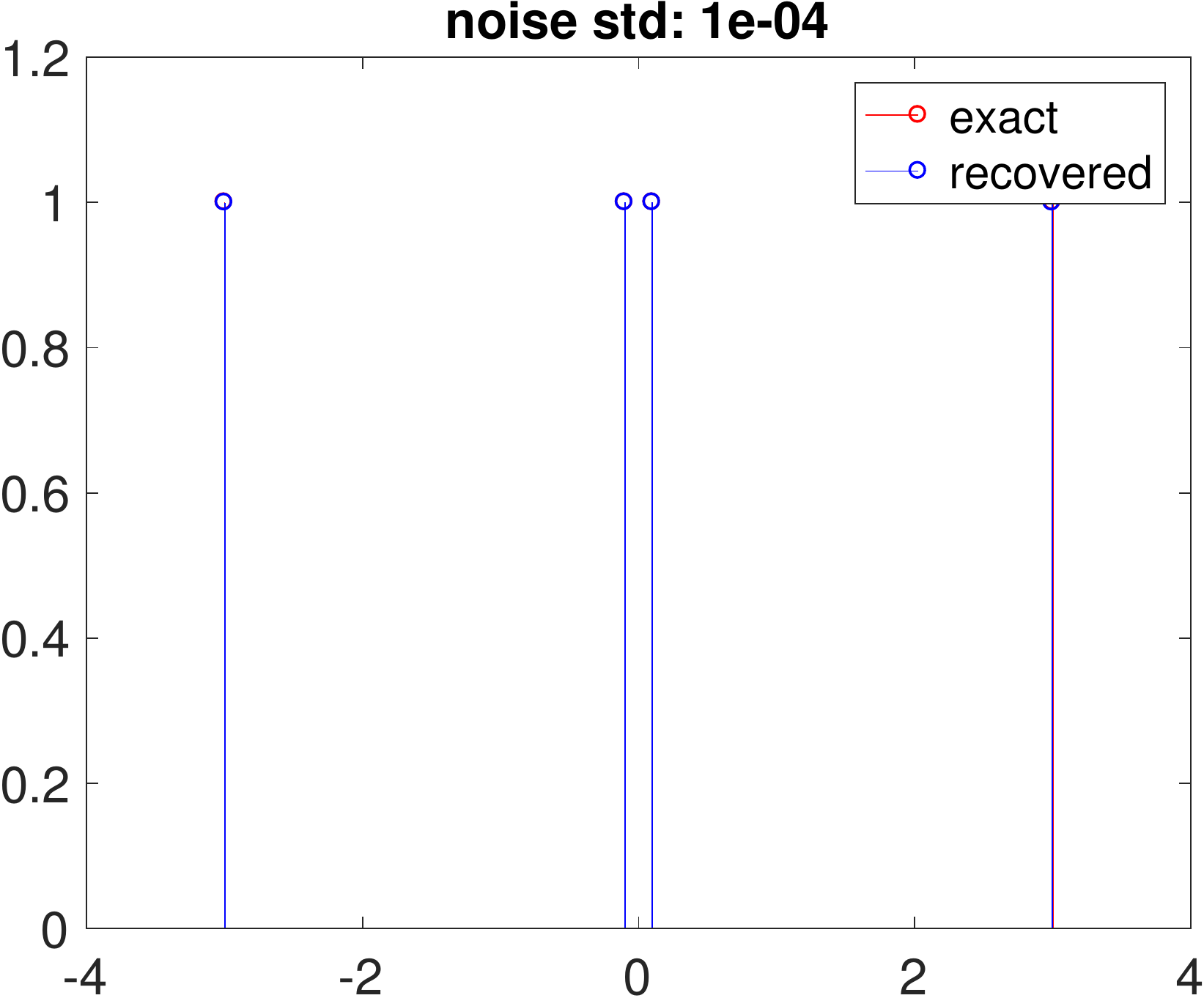} & \includegraphics[scale=0.3]{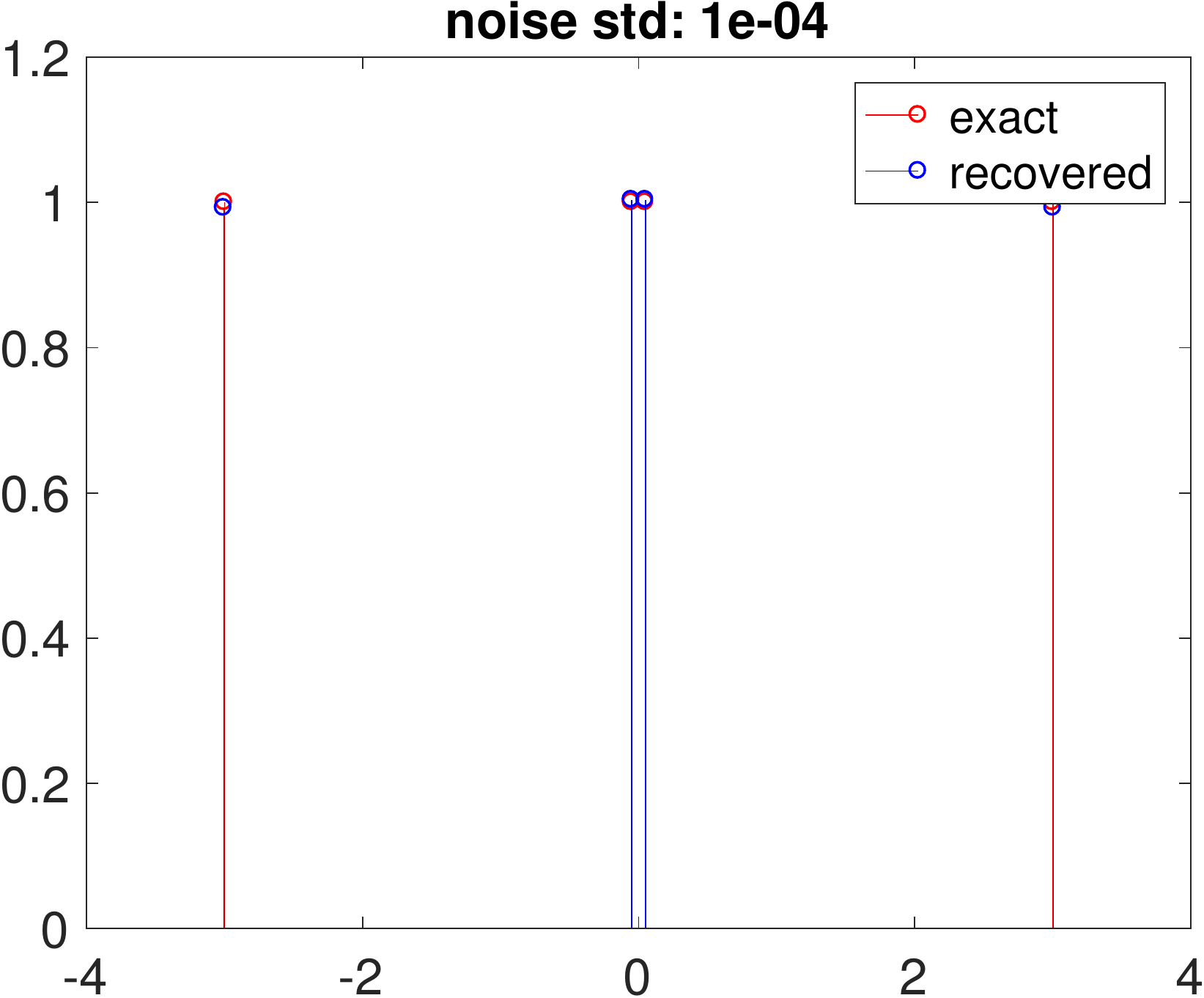}\\
    \includegraphics[scale=0.3]{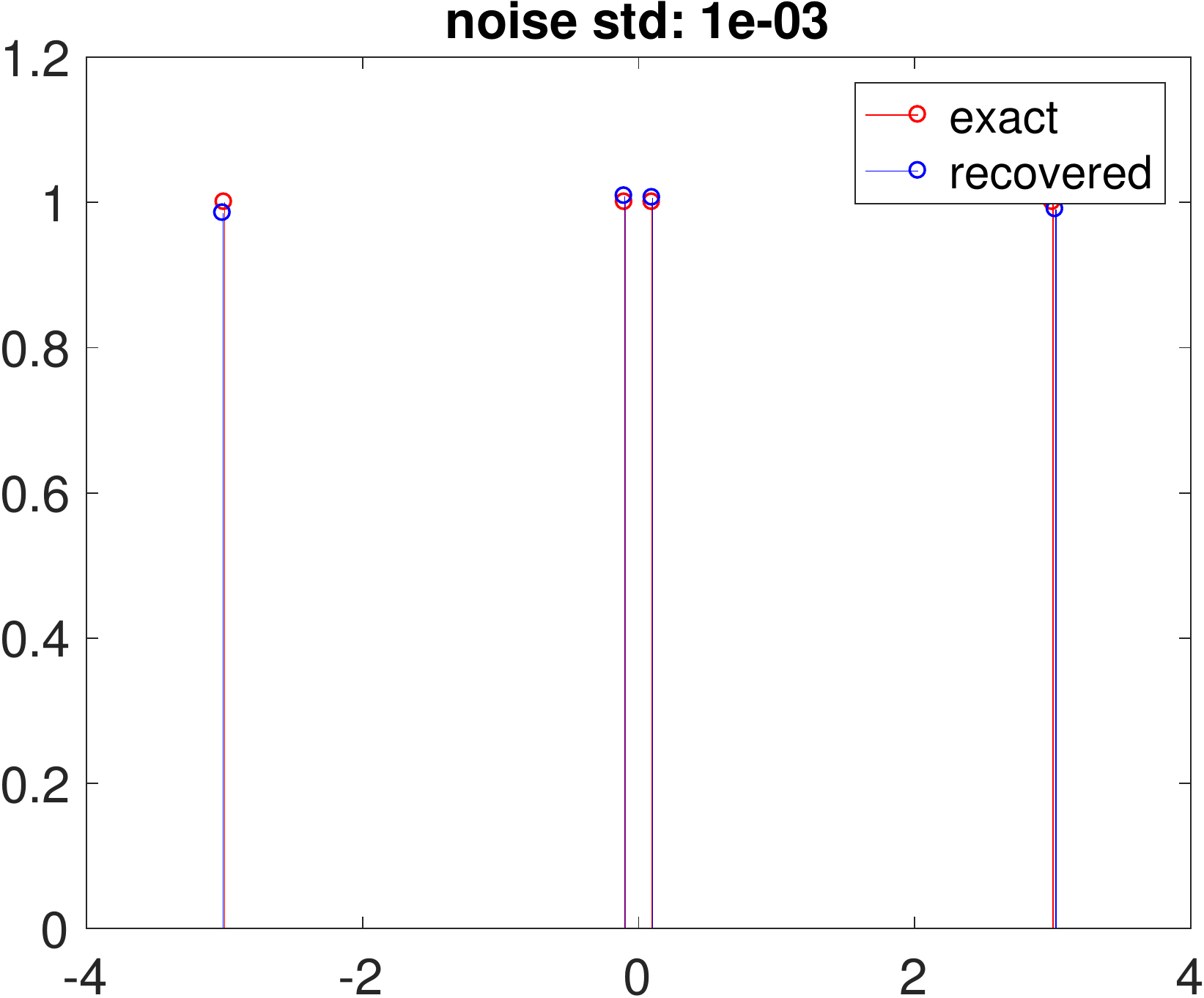} & \includegraphics[scale=0.3]{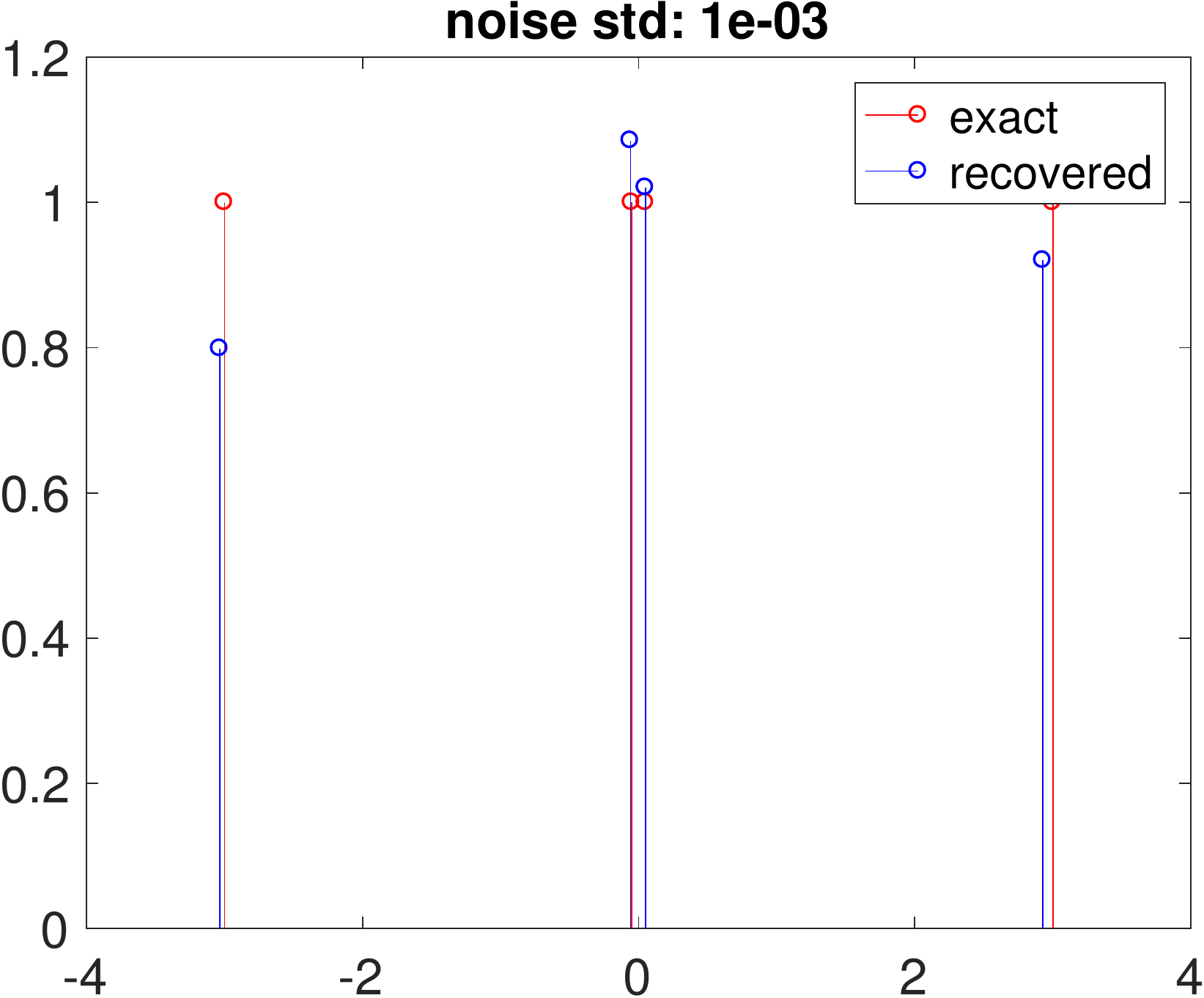}\\
    \includegraphics[scale=0.3]{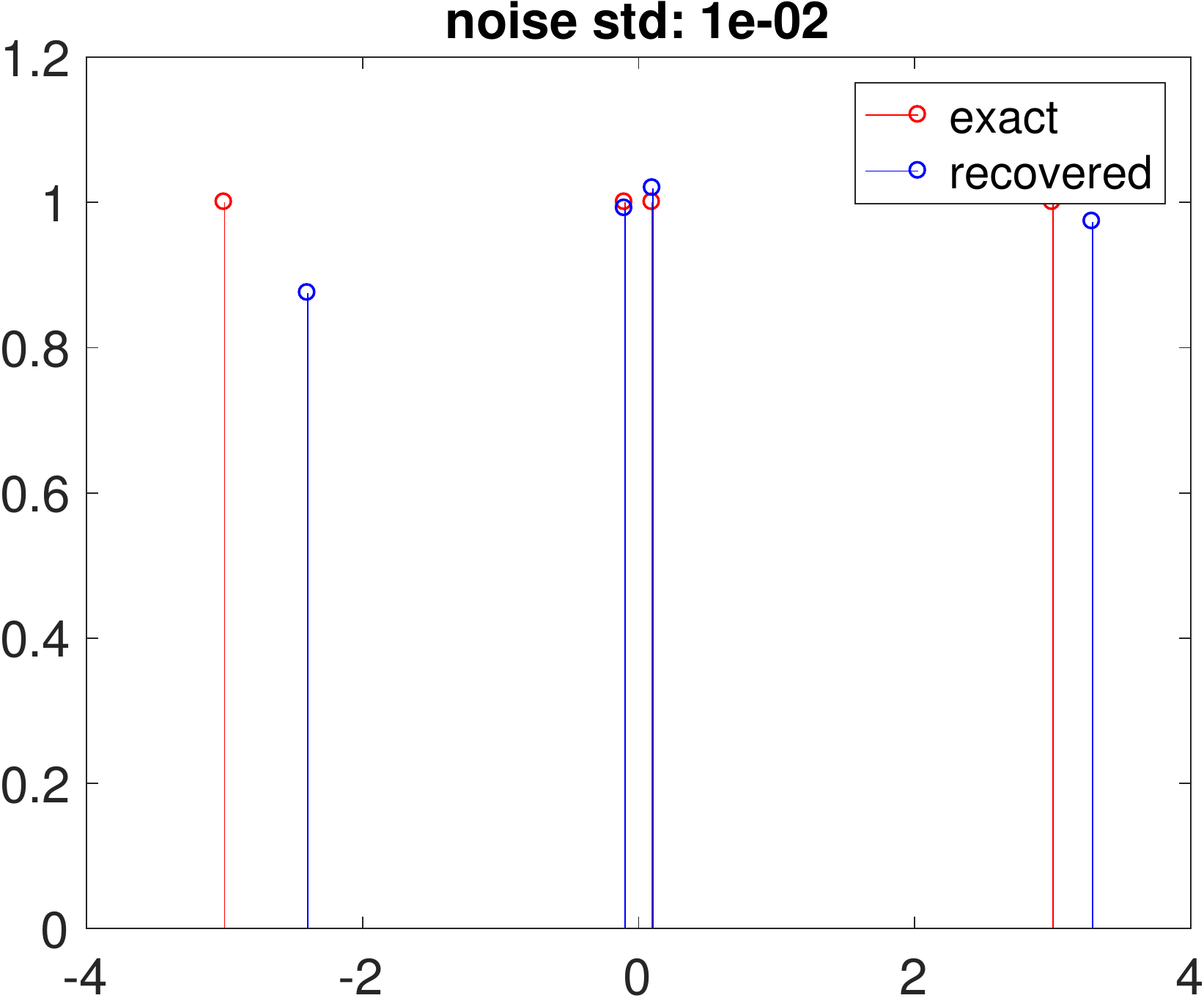} & \includegraphics[scale=0.3]{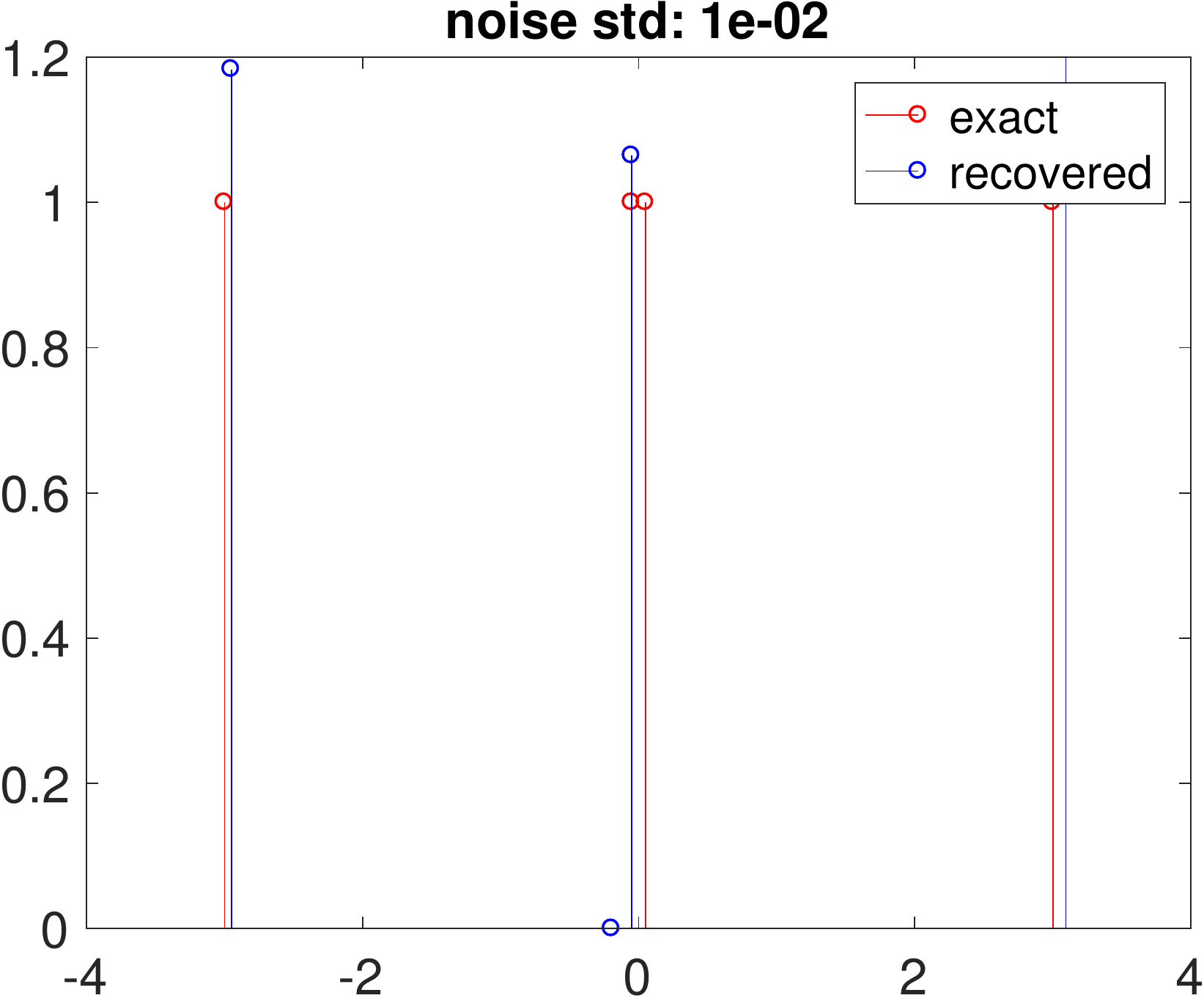}\\
    $\eps=0.1$ & $\eps=0.05$
  \end{tabular}
  \caption{Molecule case for different noise levels and different gaps.}
  \label{fig:mol}
\end{figure}

We present two examples. The inverse temperature is $\beta=100$ and the number of Matsubara points
is $N=128$. In these two examples, the gap is $\eps=0.1$ and $0.05$, respectively. For fixed values
of $\beta$ and $N$, the smaller the gap the harder the problem. The noise in $G(z_n)$ is additive
\[
G(z_n) \leftarrow G(z_n) + \sigma\cdot M \cdot N_{\C}(0,1),
\]
where $M = \left(\sum_n |G(z_n)|^2/N \right)^{1/2}$ is the average magnitude and $N_{\C}(0,1)$ is
the standard complex normal distribution. The chosen noise levels are $\sigma=10^{-4}$, $10^{-3}$,
and $10^{-2}$.  The results are summarized in Figure \ref{fig:mol}, where the two columns correspond
to $\eps=0.1$ and $0.05$, respectively.
\begin{itemize}
\item At $\sigma=10^{-4}$, the algorithm gives perfect reconstruction for both gaps.
\item At $\sigma=10^{-3}$, the algorithm provides perfect reconstruction for $\eps=0.1$, while
  some error for $\eps=0.05$.
\item At $\sigma=10^{-2}$, there are significant errors for both gap values.
\end{itemize}

\section{Condensed matter case}\label{sec:cdm}

Recall that in the condensed matter case $A(x)>0$ and the prior is that the spectral function can be
well-approximated by a number of quasi-particles.

The first step is to construct an accurate interpolant of $G(z)$ for $z\in [ai,bi]$ where
$a=\frac{\pi}{\beta}$ and $b=\frac{(2N-1)\pi}{\beta}$. To motivate the interpolant, consider the
case of a quasi-particle at location $-pi$ for $p>0$, i.e., $G(z) \approx \frac{1}{z+pi}$. When $p$
is close to zero, $G(z)$ becomes quite steep when $z \in i\R^+$ approaches the origin, making
interpolation difficult. The key idea is to consider $1/G(z) \approx z+pi$, which is easy for
interpolation.

Let $H(z)=1/G(z)$. At the Matsubara points $\{z_n\}$, we hold $H(z_n)= 1/G(z_n)$. By constructing a
high-order spline interpolant $\tH(z) \approx H(z)$ in $[ai,bi]$ based on the data
$\{(z_n,H(z_n)\}$, the interpolant $\tG(z)$ of $\G(z)$ in $[ai,bi]$ is defined as
\begin{equation}
  \tG(z) \equiv 1/\tH(z).
  \label{eq:cdmint}
\end{equation}

\begin{figure}[h!]
  \begin{tabular}{ccc}
    \includegraphics[scale=0.28]{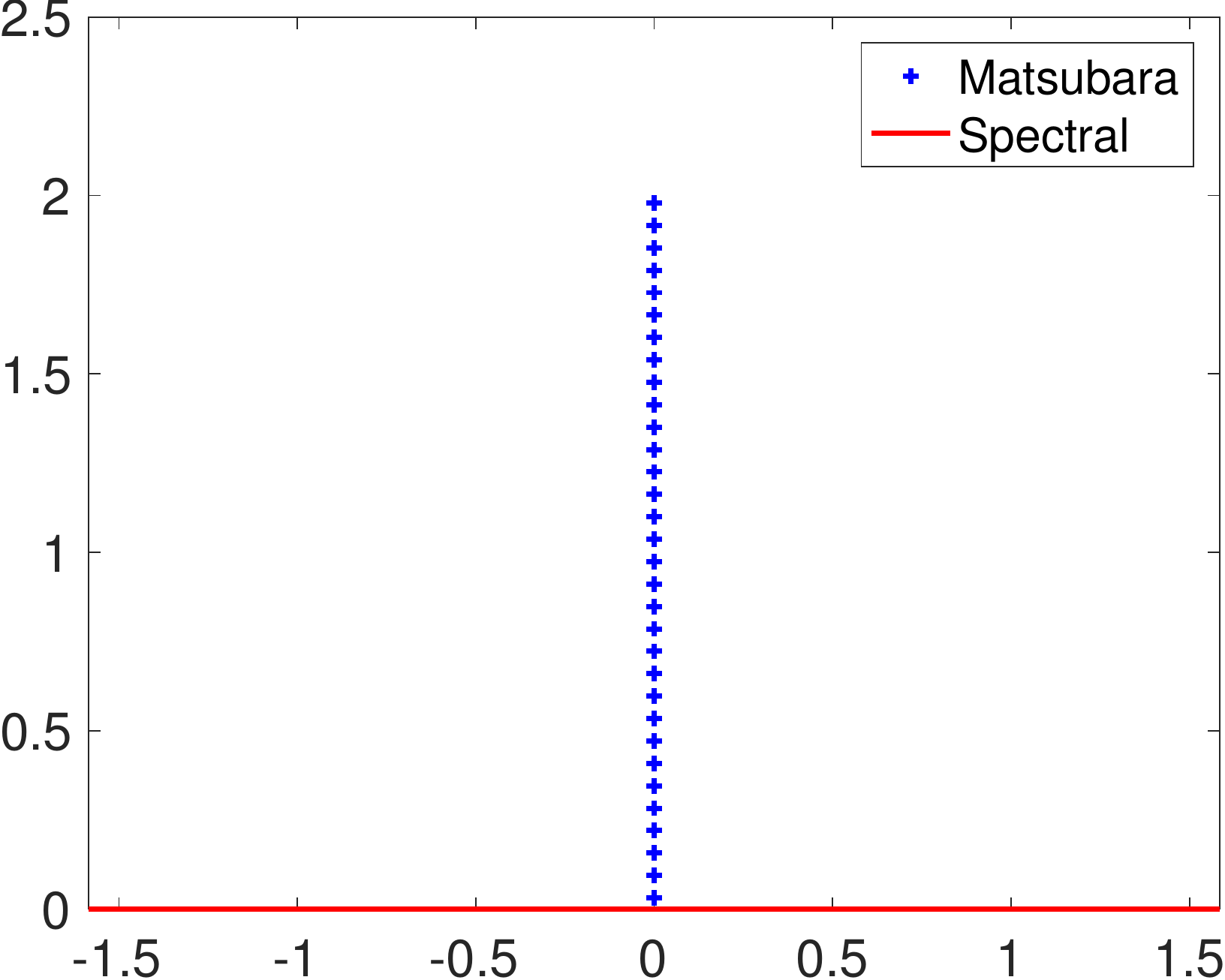}&
    \includegraphics[scale=0.28]{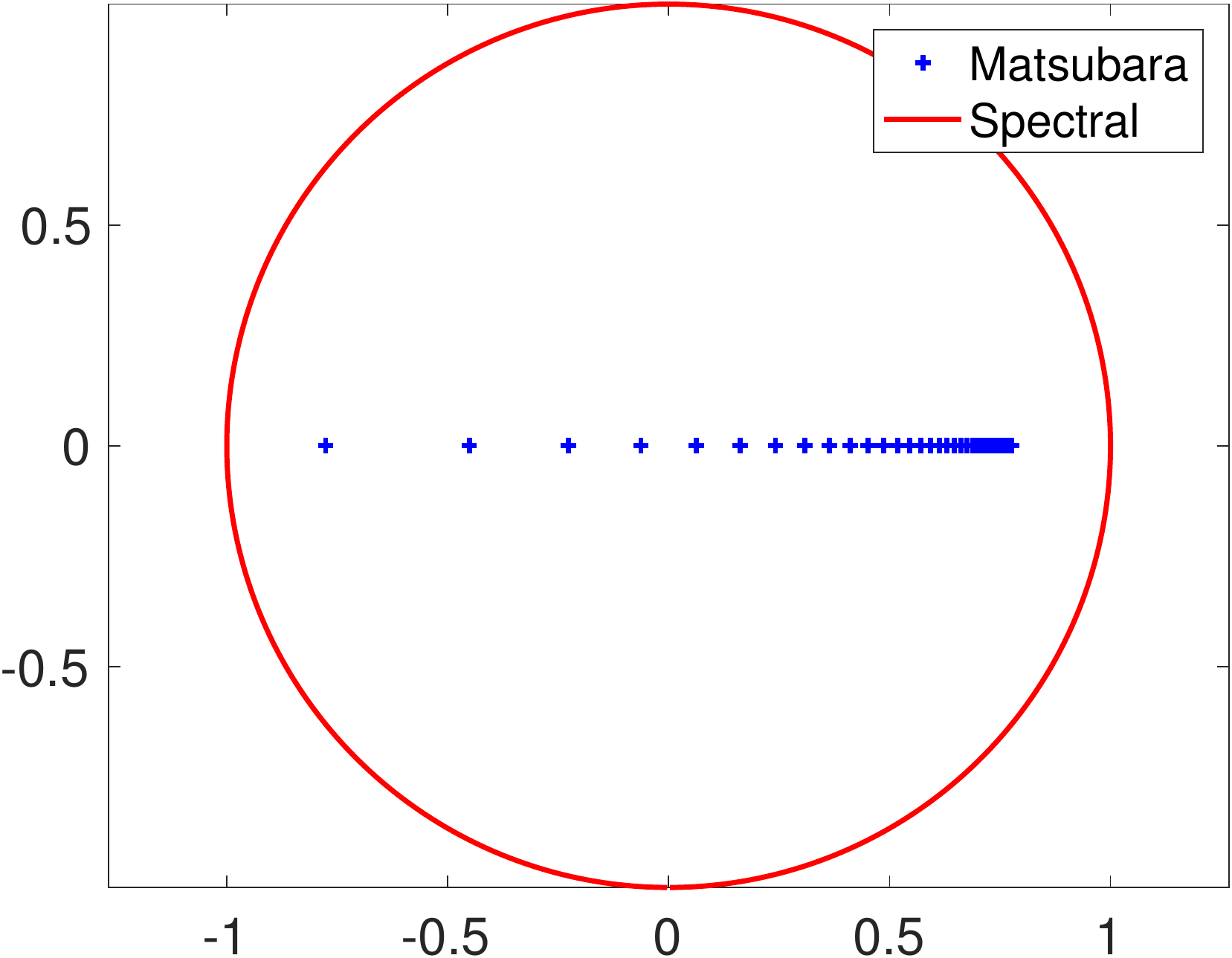}&
    \includegraphics[scale=0.28]{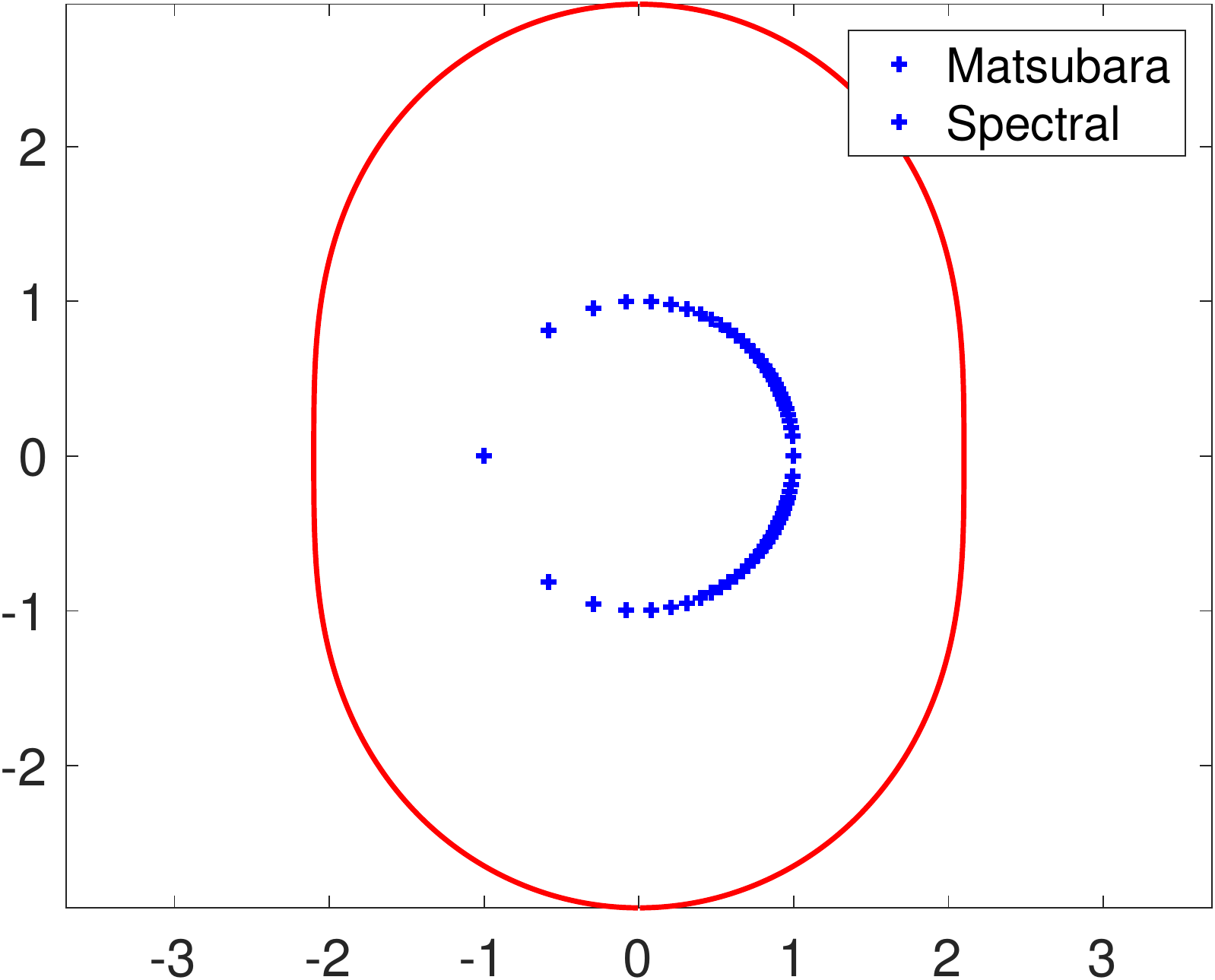}\\
    $z$ & $w$ & $t$
  \end{tabular}
  \caption{Conformal maps from $z$ to $w$ and to $t$, unzipping the interval $[ai,bi]$ into a unit
    disk.}
  \label{fig:cdmunzip}
\end{figure}

The second step is to identify the poles of the quasi-particles. Let $q=\sqrt{ab}$ and introduce the
following sequence of conformal mappings from $z$ to $w$ and to $t$ that unzip the interval
$[ai,bi]$ in the $z$ plane to the unit disk $\D$ in the $t$ plane (see Figure \ref{fig:cdmunzip})
\begin{equation}
  w = \frac{z-qi}{z+qi},\quad  t = \frac{w}{r} + \sqrt{ \frac{w^2}{r^2} -1 },
  \label{eq:cdmmap}
\end{equation}
where $r=(b-q)/(b+q)$. The inverse maps are
\begin{equation}
  z = -qi \frac{w+1}{w-1}, \quad w = \frac{r}{2} \left(t+\frac{1}{t}\right).
  \label{eq:cdmmapinv}
\end{equation}
In the $t$ plane, the function $G(t)\equiv G(z(w(t)))$ is analytic outside $\D$ and takes the form
\[
G(t) = \sum_j \frac{T_j}{t-\t_j} + f(t),
\]
where $\{\t_j\}$ are the poles outside $\D$ and $f(t)$ is analytic outside $\D$. Since the poles of
$G(z)$ (i.e., the locations of the quasi-particles) map to the poles in $G(t)$ outside $\D$, it is
sufficient to find $\{\t_j\}$ in the $t$ plane outside $\D$. Following the discussion in Section
\ref{sec:mol}, Prony's method identifies the poles $\{\t_j\}$ outside $\D$. Applying the inverse
maps \eqref{eq:cdmmapinv} from $t$ to $w$ to $z$ gives the locations of the quasi-particles
$\{\z_j\}$ in the $z$ plane.

In the third step, we solve the constrained optimization problem
\begin{equation}
  \min_{A_j: \forall x\;\Im(\sum_j A_j/(x-\z_j))\le 0} \sum_i \left|G(z_i) - \frac{1}{2\pi}\sum_j \frac{A_j}{z_i-\z_j} \right|^2
  \label{eq:cdmopt}
\end{equation}
to compute the weights $\{A_j\}$ of the quasi-particles. Finally, $A(x)$ can be approximated by
evaluating
\begin{equation}  
  -2 \Im \frac{1}{2\pi} \sum_j \frac{A_j}{x+i\eta - \z_j}
  \label{eq:cdmfnl}
\end{equation}
for a sufficiently small positive $\eta$. The constraint in \eqref{eq:cdmopt} is included to ensure
that \eqref{eq:cdmfnl} is positive.

To implement this algorithm, we need to take care several numerical issues.
\begin{itemize}
\item For the spline interpolation for $H(z)$, a 5th order spline is used.

\item To compute the Fourier coefficients $\{\hG_k\}$, we again use a uniform grid
  $\theta_n=\frac{2\pi n}{N_s}$ for $n=0,\ldots,N_s-1$.  $N_S$ is chosen such that
  \[
  h \ll \sqrt{\frac{a}{b}}, \quad\text{i.e.,}\quad
  N_s \gg \sqrt{\frac{b}{a}}.
  \]
  to ensure the exponential convergence of the trapezoidal rule.

\item The constrained optimization problem \eqref{eq:cdmopt} is solved with CVX \cite{cvx}.

\end{itemize}

\subsection{Numerical results}

\begin{figure}[h!]
  \begin{tabular}{cc}
    \includegraphics[scale=0.3]{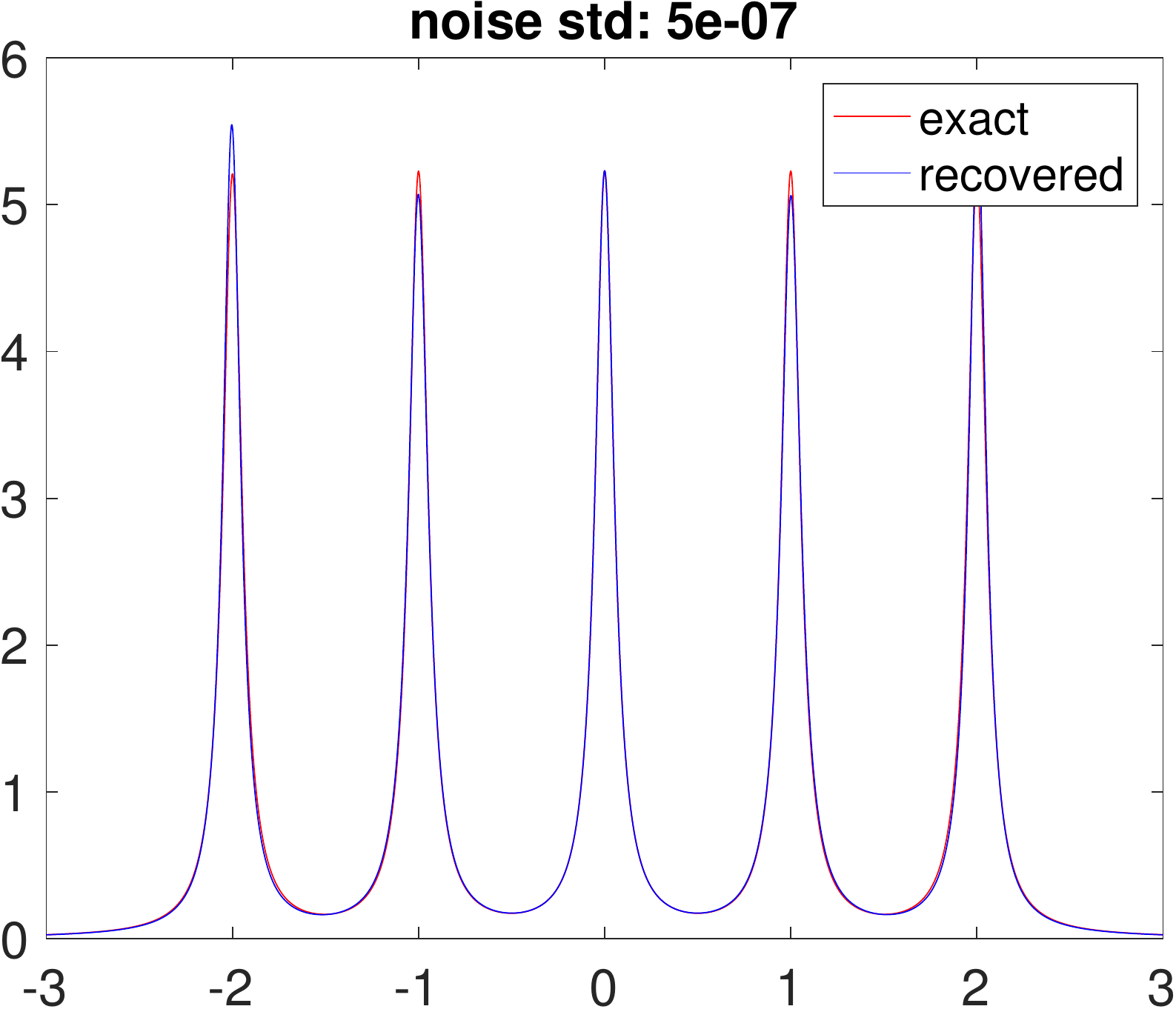} & \includegraphics[scale=0.3]{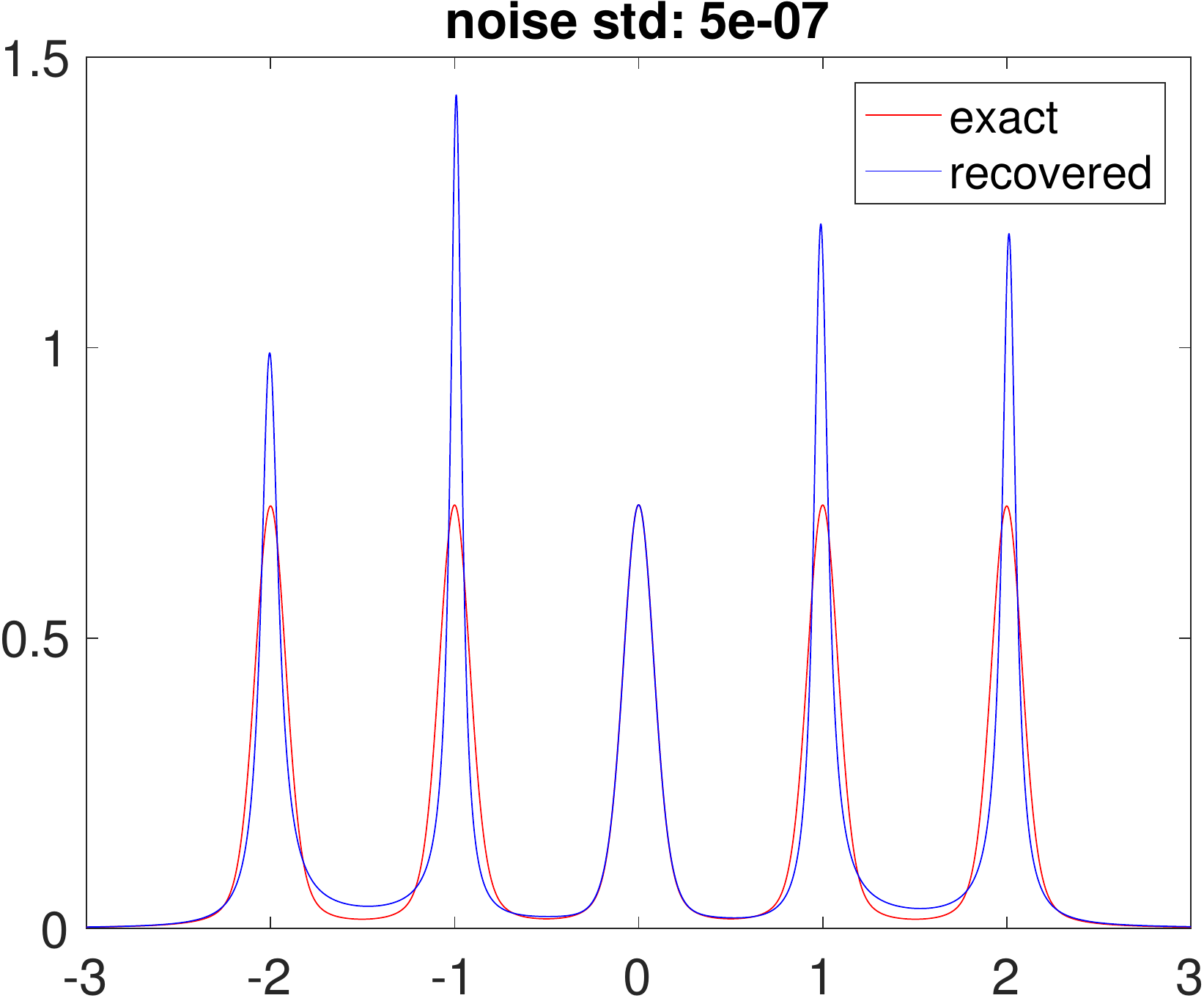}\\
    \includegraphics[scale=0.3]{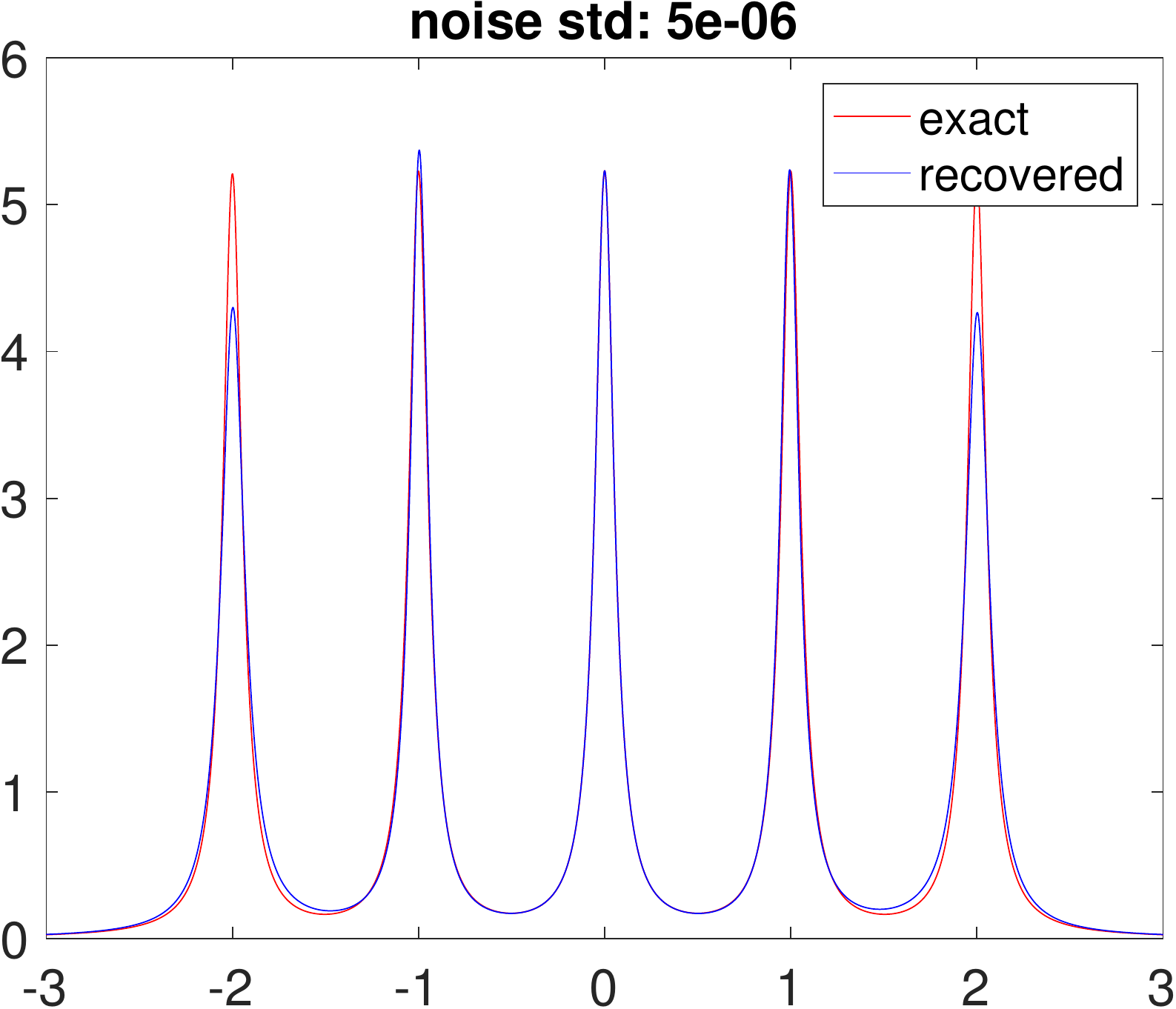} & \includegraphics[scale=0.3]{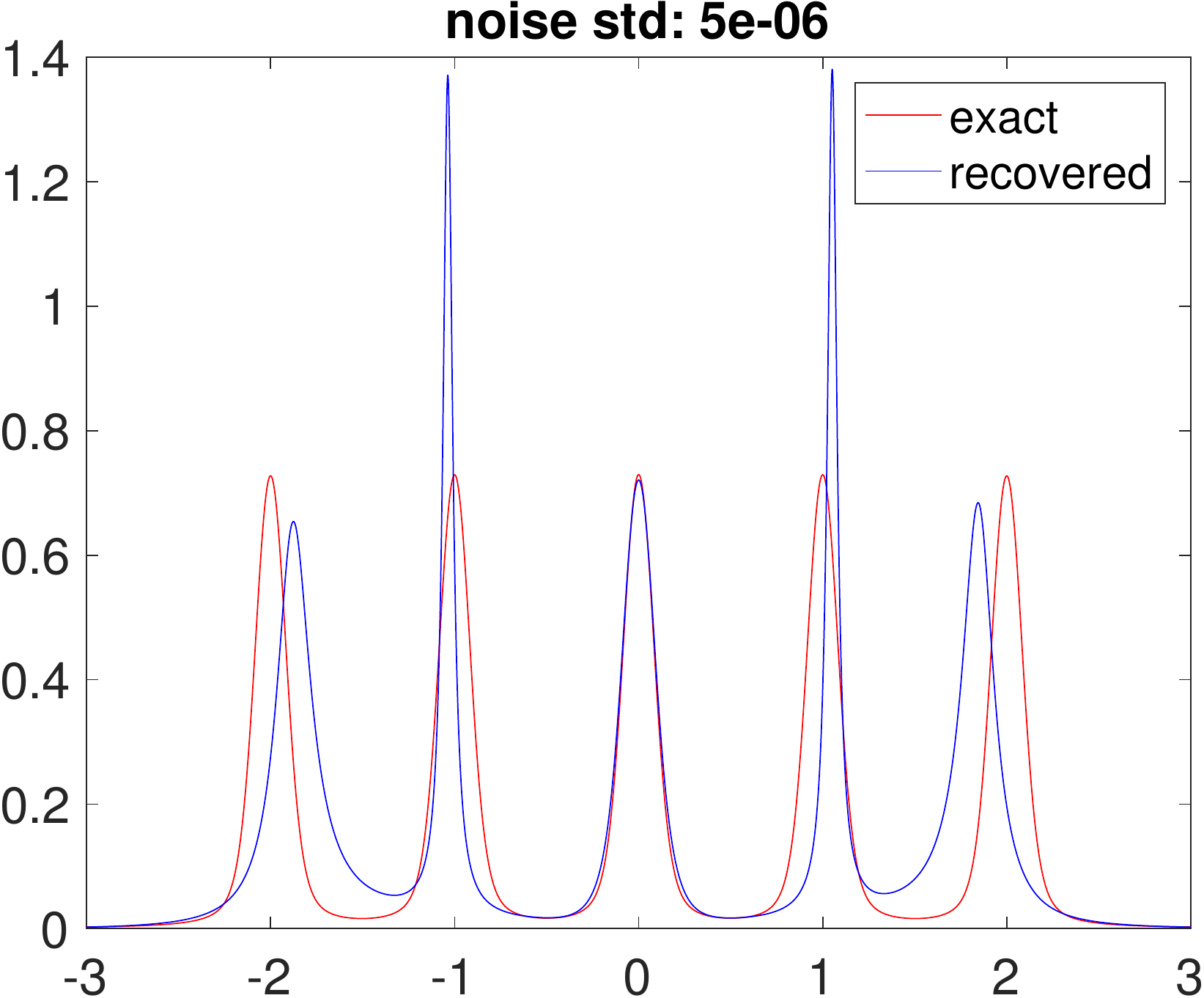}\\
    \includegraphics[scale=0.3]{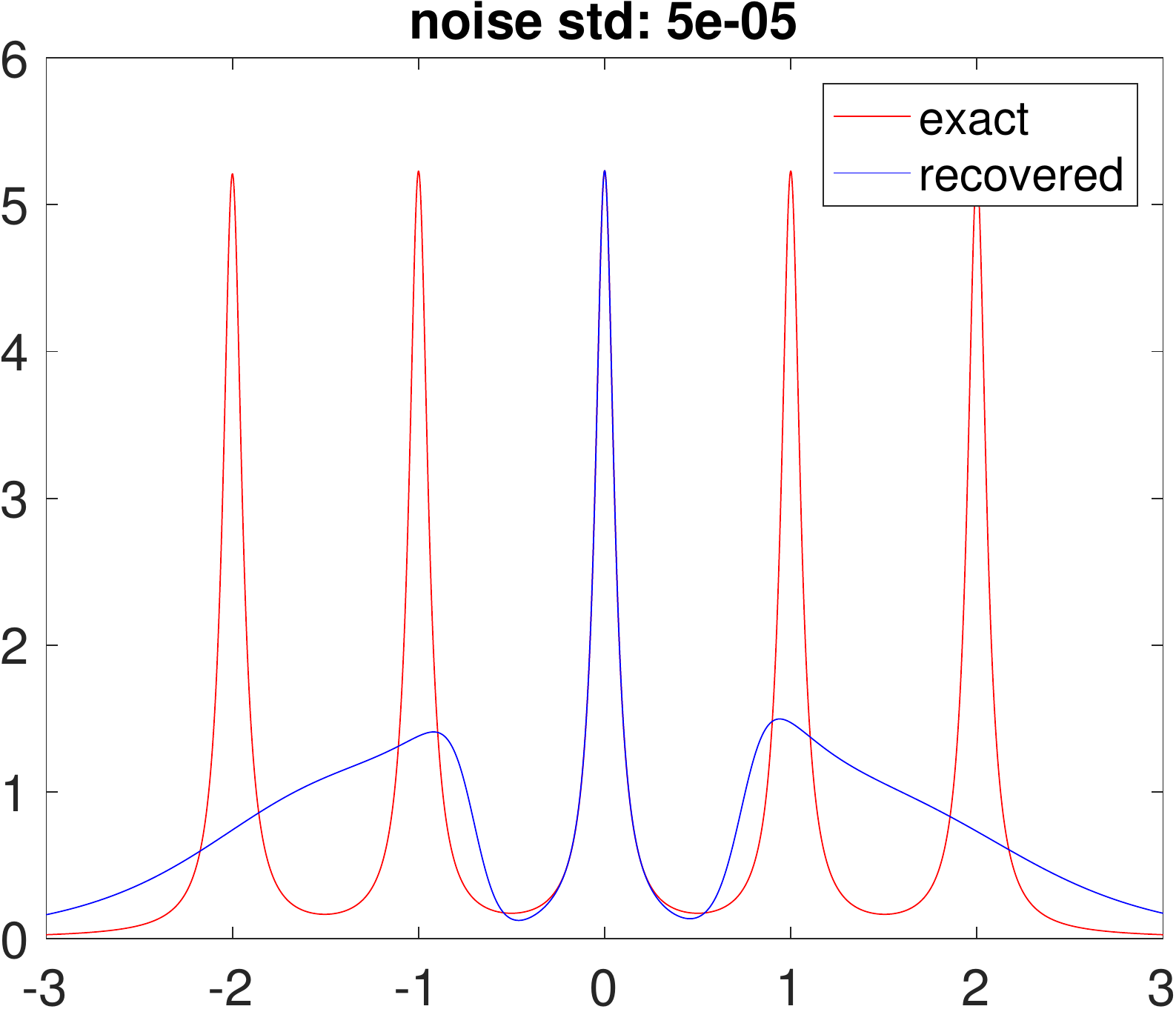} & \includegraphics[scale=0.3]{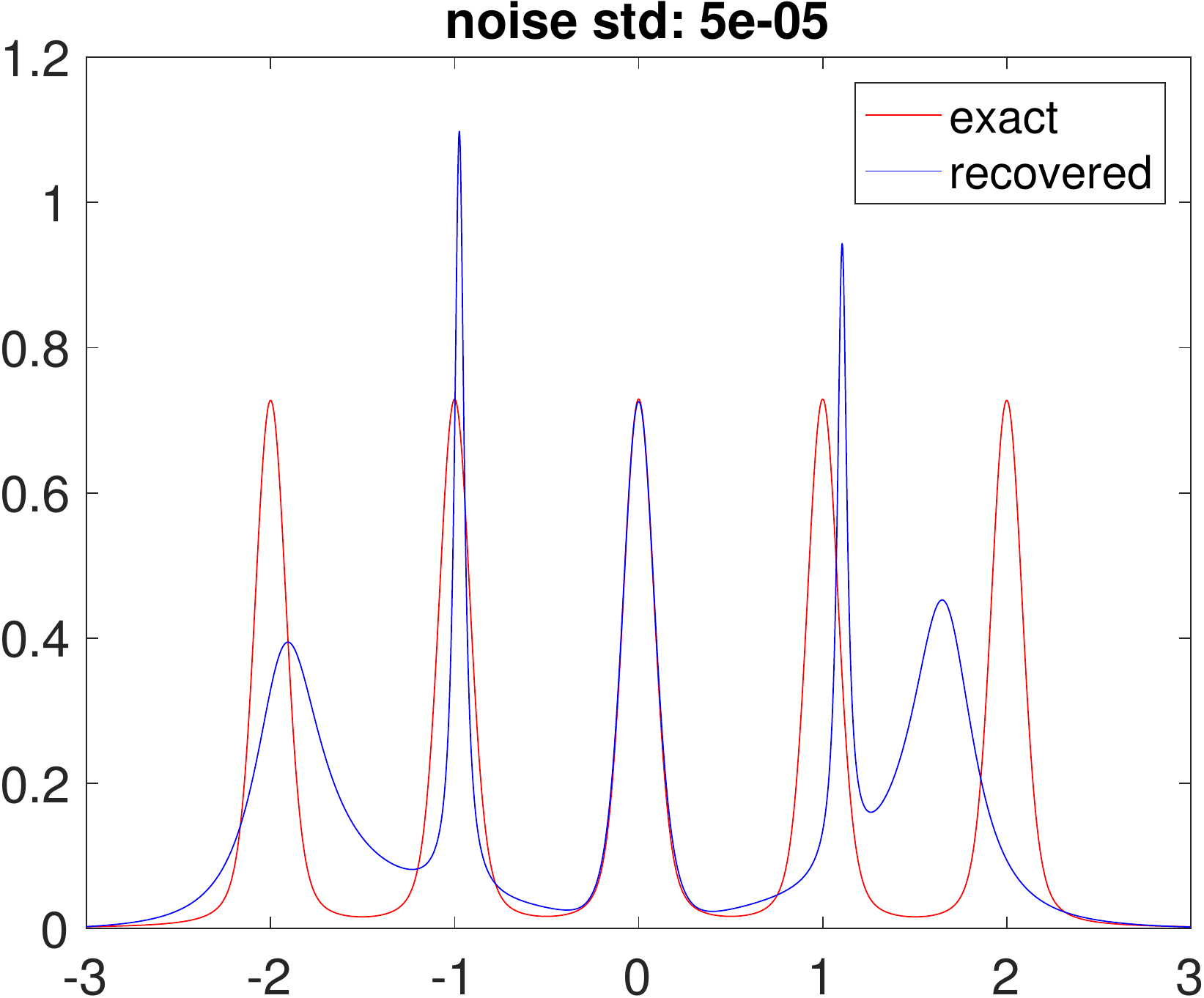}\\
    Quasi-particles & Gaussians
  \end{tabular}
  \caption{Condensed matter case for different noise levels and different spectral distribution models.}
  \label{fig:cdm}
\end{figure}

We present two numerical examples. The inverse temperature $\beta=100$ and the number of Matsubara
points $N=256$. In the first example, $A(x)$ indeed corresponds to a sum of quasi-particles at
\[
\{-2-0.03i,-1-0.03i,0-0.03i,1-0.03i,2-0.03i\} 
\]
and here the quasi-particle prior used by the algorithm is correct. In the second example, $A(x)$ is
a sum of five Gaussians centered at
\[
\{-2,-1,0,1,2\}
\]
with variance equal to $\frac{1}{200}$ and here the prior is thus misspecified. The noise in
$G(z_n)$ is again
\[
G(z_n) \leftarrow G(z_n) + \sigma\cdot M \cdot N_{\C}(0,1),
\]
where $M = \left(\sum_n |G(z_n)|^2/N \right)^{1/2}$ is the average magnitude and $N_{\C}(0,1)$ is
the standard complex normal distribution. The chosen noise levels are $\sigma=5\cdot 10^{-7}$,
$5\cdot 10^{-6}$, $5\cdot 10^{-5}$. The results are summarized in Figure \ref{fig:cdm}, where we
plot $-2 \Im G(x+i\eta)$ on a horizontal line ($\eta=0.01$) close to the real axis.

\begin{itemize}
\item At $\sigma=5\cdot 10^{-7}$, the algorithms gives a perfect reconstruction for the
  quasi-particle example. For the Gaussian example, the peak locations are well identified but the
  widths and heights are a bit off.
\item At $\sigma=5\cdot 10^{-6}$, the quasi-particle example still shows a good reconstruction. For
  the Gaussian example, there is a shift (towards the center) for the Gaussians away from the
  origin.
\item At $\sigma=5\cdot 10^{-5}$, the noise level is too large for the algorithm in both examples.
\end{itemize}

\bibliographystyle{abbrv}

\bibliography{ref}

\end{document}